\documentclass[11pt,reqno]{amsart}
\usepackage{amsmath,amssymb,amsxtra,latexsym,amsthm,amscd,amsfonts,mathrsfs,pb-diagram}
\usepackage[colorlinks=true,linkcolor=magenta,citecolor=blue]{hyperref}
\usepackage[margin=1in]{geometry}
\usepackage{enumitem}
\usepackage{mathabx}
\usepackage{epsfig}

\usepackage{color}

\begin{document}
	\title[\hfil evolution equations with modified Hartree Nonlinearity]{A note on evolution equations with modified Hartree Nonlinearity}
\thanks{}


	\maketitle
	\centerline{\scshape Khaldi Said}
	\medskip
	
	{\footnotesize
	\centerline{Laboratory of Analysis and Control of PDEs, Djillali Liabes University}
	\centerline{P.O. Box 89, Sidi-Bel-Abbes 22000, Algeria}
	\centerline{Email: saidookhaldi@gmail.com}
	}
\medskip

\centerline{\scshape}
\medskip
{\footnotesize
	\centerline{} 
	\centerline{}
	\centerline{}
} 
% Do not forget to end the {\footnotesize by the sign }

\begin{abstract}
We introduce a mathematical model in $\mathbb{R}^{n}$ for evolution equations with modified generalized Hartree nonlinearity given by $S_{\alpha,p,q}(u)=I_{\alpha}(|u|^{p+q}).$ 
One can see that this nonlinearity is not integrable due to the boundedness property of Riesz potential. In other words, we cannot deal with the Cauchy problem of semi-linear evolution equations with $S_{\alpha,p,q}(u)$ and $L^{1}$-initial velocity.

We will show that $S_{\alpha,p,q}(u)$ produces the same semi-critical exponent that guarantees the global existence of small data solutions as in the well known generalized Hartree nonlinearity $H_{\alpha,p,q}(u)=|u|^{p}I_{\alpha}(|u|^{q})$ provided that the initial velocity belongs to $L^{m}(\mathbb{R}^{n})$, with $m>1$.

We can expect a relation between some physical systems (Schr\"{o}dinger–Poisson) that are modeled and solved using Hartree nonlinearity and those in their modified form due to this coincidence property in the semi-critical exponent.
\end{abstract}
\numberwithin{equation}{section}
\newtheorem{theorem}{Theorem}[section]
\newtheorem{lemma}[theorem]{Lemma}
\newtheorem{definition}[theorem]{Definition}
\newtheorem{remark}[theorem]{Remark}
\newtheorem{proposition}[theorem]{Proposition}
\allowdisplaybreaks	
\section{Introduction}\label{intro}

In this research we will look at the link between two types of nonlinearities from the point of view of the equivalence of semi-critical exponent that guarantees the global existence of small data solutions to the associated evolution equation. To do this, let us consider the following semi-linear evolution model
\begin{equation}\label{1}
\left\lbrace 
\begin{array}{lll}
u_{tt}(t,x)+(-\Delta)^{\sigma}u(t,x)+(-\Delta)^{\sigma/2}u_{t}(t,x)= N(u)(t,x) \ \ \ \hfill& x\in \mathbb{R}^{n}, \ t>0,
&\cr
\\
u(0,x)=u_{0}(x), \ \ u_{t}(0,x)=u_{1}(x) & x\in \mathbb{R}^{n}.
\end{array}\right.  
\end{equation}
The linear part in model (\ref{1}) when $N(u)\equiv 0$ and $\sigma\geq 1$ is called $\sigma$-evolution equations with critical structural damping, and when $\sigma=1$, it is a linear structurally damped waves.

The nonlinearity $N(u)$ takes the following two representations, the generalized Hartree nonlinearity
\begin{align}\label{2}
N(u)(t,x) &=H_{\alpha,p,q}(u)(t,x)\nonumber\\
&=|u(t,x)|^{p}I_{\alpha}(|u|^{q})(t,x) \nonumber \\
&=|u(t,x)|^{p}\int_{\mathbb{R}^{n}}^{}\frac{|u(t,y)|^{q}}{|x-y|^{n-\alpha}}dy=\int_{\mathbb{R}^{n}}^{}\frac{|u(t,x)|^{p}|u(t,y)|^{q}}{|x-y|^{n-\alpha}}dy,
\end{align}
and its modified form
\begin{align}\label{3}
N(u)(t,x) &=S_{\alpha,p,q}(u)(t,x)\nonumber\\
&=I_{\alpha}(|u|^{p+q})(t,x)=\int_{\mathbb{R}^{n}}^{}\frac{|u(t,y)|^{p+q}}{|x-y|^{n-\alpha}}dy,
\end{align} 
with the same parameters $p,q \geq 1, \ \alpha\in (0,n).$
The nonlinearity $H_{2,1,2}(u)$ is well known in the literature as original Hartree type nonlinearity in $\mathbb{R}^{3}$ that was first written by D.R. Hartree in his research on the wave mechanics of an atom, see his papers [13, 14, 15] in \cite{dabbicco hartree} and references therein. The nonlinearity $S_{\alpha, p,q}(u)$ is first introduced here to show its close relation with the generalized Hartree nonlinearity. Of course, if we allow $\alpha \to 0$, then both (\ref{2}), (\ref{3}) become the usual power nonlinearity $ |u(t,x)|^{p+q}$.

Let us now review some previously published results for the critical exponent of a model (\ref{1}) with usual power nonlinearity  
\begin{equation}\label{4}
\left\lbrace 
\begin{array}{lll}
u_{tt}(t,x)+(-\Delta)^{\sigma}u(t,x)+(-\Delta)^{\sigma/2}u_{t}(t,x)= |u(t,x)|^{p} \ \ \ \hfill& x\in \mathbb{R}^{n}, \ t>0,
&\cr
\\
u(0,x)=u_{0}(x), \ \ u_{t}(0,x)=u_{1}(x) & x\in \mathbb{R}^{n}
\end{array}\right.  
\end{equation}
where $p>1$ and $(-\Delta)^{\sigma}$ denotes the fractional Laplacian operator with symbol $|\xi|^{2\sigma}$, i.e., by the Fourier transform $\mathcal{F}$ we have
$$
\mathcal{F}\left( (-\Delta)^{\sigma}f\right) = |\xi|^{2\sigma}\mathcal{F}\left(f\right) (\xi), \ \  \xi\in\mathbb{R}^{n}, \ \  |\xi|=(\xi_{1}^{2}+\cdots+\xi_{n}^{2})^{1/2}.$$

For example, in \cite{pham reissig kainane} (see also [6] in \cite{dabbicco hartree}), the authors found the following critical exponent $p_{crit}$ of the Cauchy problem (\ref{4}) with additional $L^m \cap L^{r}$ regularity of initial data, where $1\leq m < r\leq \infty$
\begin{equation}\label{5}
p_{crit}(n, m,\sigma)=1+\frac{2\sigma}{\frac{n}{m}-\sigma}.
\end{equation}

This exponent is critical in the sense that it guarantees both the global (in time) existence of small data Sobolev solutions for $p > p_{crit}(n, m, \sigma)$, and blow-up in finite time for $1 < p <p_{crit}(n, m, \sigma)$.
An interesting result shows that if $m=1$, $p_{crit}(n, 1, \sigma)$ belongs to the blow up region, however if $m>1$, $p_{crit}(n, m, \sigma)$ belongs to the global existence region. If we are successful in proving one partial result, either global existence or blow up, we shall refer to $p_{crit}$ as the semi-critical exponent in this study. 

For Cauchy problem (\ref{1}) with nonlinearity (\ref{2}) and additional $L^{m}$ regularity of initial data, the author in \cite{dabbicco hartree} derived the following semi-critical exponent 
\begin{equation}\label{6}
p+q\geq 1+\frac{2\sigma+\alpha}{\frac{n}{m}-\sigma} \ \ \ \ if \ m>1, \  \ p+q>1+\frac{2\sigma+\alpha}{n-\sigma} \ \ if \ m=1,
\end{equation}
this only guarantees the global (in time) existence of small data solutions. It is still unclear whether the counterpart of condition (\ref{6}) causes the solution to blow up or not due to the nonlocal (in space) property, see the end of page 15 in reference [6] cited in \cite{dabbicco hartree}.

Fortunately, there are two approaches for determining the semi-critical exponents of some semi-linear evolution equations. The first is to establish the global existence result using certain sharp linear estimates and the Banach fixed point theorem, while the second is to prove the blow-up result using a test function based on contradiction argument.

Without going into more details, in this research we will rely on the first way to prove not only the global (in time) existence of small data solutions to (\ref{1}), but also to show the close relation between nonlinearities (\ref{2}) and (\ref{3}).

Finally, the idea of researching this type of coincidence comes from the fact that if someone successfully investigates two different evolutionary equations and then discovers the same result exists under the same parameters, he can conjecture that they have a hidden relation. We can cite a powerful example the so-called diffusion phenomenon between classical damped waves and heat equations.

\section{Useful tools}\label{Main tools}
First,  a constant $c>0$ such that $f \leq cg$ does not plays any role in our analysis, we omit it and we write $f\lesssim g$. We now introduce the Riesz potential and the boundedness property from $L^{q}(\mathbb{R}^{n})$ to $L^{r}(\mathbb{R}^{n})$ spaces known as Hardy-Littlewood-Sobolev inequality. For more details see the pioneering book \cite[Page 119]{stein}.
\begin{definition} Let $\alpha \in(0,n)$. We formally define the normalized Riesz potential as:
	$$(I_{\alpha}f)(x):=\frac{\Gamma((n-\alpha)/2)}{\pi^{n/2}2^{\alpha}\Gamma(\alpha/2)} \int_{\mathbb{R}^{n}}^{}\frac{f(y)}{|y-x|^{n-\alpha}}dy, \ \ x\in \mathbb{R}^{n},$$
	where $\Gamma$ is the Euler Gamma function.
\end{definition}
\begin{lemma} \label{hardy sobolev inequality}
	If $f \in L^{q}(\mathbb{R}^{n})$ for some $q\in(1,n/\alpha)$, then  $I_{\alpha}f \in L^{r}(\mathbb{R}^{n})$ for some $q<r<\infty$ and satisfies the inequality
	$$\|I_{\alpha}f\|_{L^{r}(\mathbb{R}^{n})}\lesssim \|f\|_{L^{q}(\mathbb{R}^{n})}, \ \ \frac{1}{q}-\frac{1}{r}=\frac{\alpha}{n}.$$
	Of course, $q<r$ due to the fact that $\alpha>0$.	
\end{lemma}

We now recall the following sharp long-time decay estimates which is very important tool to demonstrate Theorem \ref{GlobalExistence1}. Here and in the following we choose $u_0=0$ for brevity.
\begin{lemma}[\cite{pham reissig kainane}]\label{linear estimates}
	Let $m\in [1,2)$. The Sobolev solutions $u$ to the linear equation in (\ref{1}) satisfy the $(L^{m}\cap L^{2})-L^{2}$ estimates:
	\begin{align} 
	\|u(t,\cdot)\|_{L^{2}(\mathbb{R}^{n})} &\lesssim (1+t)^{1-\frac{n}{\sigma}\left(\frac{1}{m}-\frac{1}{2}\right) }\|u_{1}\|_{L^{m}(\mathbb{R}^{n})\cap L^{2}(\mathbb{R}^{n})}, \label{7} 
	\end{align}
	\begin{align} 
	\|((-\Delta)^{\sigma/2}u, u_{t})(t,\cdot)\|_{L^{2}(\mathbb{R}^{n})} &\lesssim (1+t)^{-\frac{n}{\sigma}\left(\frac{1}{m}-\frac{1}{2}\right)}\|u_{1}\|_{L^{m}(\mathbb{R}^{n})\cap L^{2}(\mathbb{R}^{n})}, \label{8} 
	\end{align}
	and the $L^{2}-L^{2}$ estimates:
	\begin{align} 
	\|u(t,\cdot)\|_{L^{2}(\mathbb{R}^{n})} &\lesssim (1+t)\|u_{1}\|_{L^{2}(\mathbb{R}^{n})}, \label{9} 
	\end{align}
	\begin{align} 
	\|((-\Delta)^{\sigma/2}u, u_{t})(t,\cdot)\|_{L^{2}(\mathbb{R}^{n})} &\lesssim \|u_{1}\|_{L^{2}(\mathbb{R}^{n})}. \label{10} 
	\end{align}
\end{lemma}
\smallskip
We also need to the fractional Gagliardo-Nirenberg inequality in the following lemma.
\begin{lemma} [\cite{pham reissig kainane}]\label{FGN}
	Let $1<q<\infty$, $\sigma>0$. Then, the following fractional Gagliardo-Nirenberg inequality holds for all $u\in H^{\sigma}(\mathbb{R}^{n})$
	\[\|u\|_{L^{q}(\mathbb{R}^{n})}\lesssim \|(-\Delta)^{\sigma/2}u\|_{L^{2}(\mathbb{R}^{n})}^{\theta_{q}}\,\|u\|_{L^{2}(\mathbb{R}^{n})}^{1-\theta_{q}},\]
	where \[\theta_{q}=\frac{n}{\sigma}\left(\frac{1}{2}-\frac{1}{q}\right)\in\left[0,1\right].\]
\end{lemma}
Finally, this integral inequality is used to deal with the Duhamel's formula. 
\begin{lemma}[\cite{ebert reissig}]\label{Integral inequality}
	Let $a, b\in\mathbb{R}$ such that $\max\{a,b\}>1$. Then, it holds
	$$\int_{0}^{t}(1+t-\tau)^{-a}(1+\tau)^{-b}d\tau\lesssim (1+t)^{-\min\{a,b\}}.$$ 
\end{lemma} 

\smallskip
\section{Main Result}\label{Main results}
Our main result is reads as follows.
\begin{theorem}\label{GlobalExistence1}
	Let $m\in (1,2)$, let us consider the Cauchy problem (\ref{1}) with nonlinearity (\ref{3}), where $\sigma\geq1$, $p,q \geq 1$ and $\alpha\in (0,n)$. We assume the following conditions for $p,q$ and the dimension $n$:
	\begin{equation}\label{11}
	\left\lbrace 
	\begin{array}{lll}
	\frac{2}{m}+\frac{2\alpha}{n}\leq p+q \leq \frac{n+2\alpha}{n-2\sigma} \ \ \ \hfill& {if } \ \ \ \ 2\sigma<n\leq \frac{2\sigma+2\sqrt{\sigma(\sigma+m(2-m)\alpha)}}{2-m},
	&\cr
	\\
	\frac{2}{m}+\frac{2\alpha}{n}\leq p+q < \infty & {if } \ \
	n \leq 2\sigma.
	\end{array}\right.  
	\end{equation}
	Moreover, we suppose
	\begin{equation}\label{12}  
	p+q>1+\frac{2\sigma+\alpha}{\frac{n}{m}-\sigma}.
	\end{equation}\\
	Then, there exists a constant $\varepsilon_0>0$ such that for any data 
	$$u_{1}\in L^{m}(\mathbb{R}^{n})\cap L^{2}(\mathbb{R}^{n}) \ \textit{with}\  \left\|u_{1}\right\|_{L^{m}(\mathbb{R}^{n})\cap L^{2}(\mathbb{R}^{n})}<\varepsilon_0,$$ 
	we have a uniquely determined globally (in time) solution
	$$u\in\mathcal{C}\left([0,\infty), H^{\sigma}(\mathbb{R}^{n})\right)\cap \mathcal{C}^{1}\left([0,\infty),L^{2}(\mathbb{R}^{n})\right)$$
	to (\ref{1}) and (\ref{3}). Furthermore, the solution satisfies the estimates:
	\begin{equation*}
	\left\|u(t,\cdot)\right\|_{L^{2}(\mathbb{R}^{n})} \lesssim (1+t)^{1-\frac{n}{\sigma}\left(\frac{1}{m}-\frac{1}{2} \right)}\left\|u_{1}\right\|_{L^{m}(\mathbb{R}^{n})\cap L^{2}(\mathbb{R}^{n})},
	\end{equation*}
	\begin{equation*}
	\left\|((-\Delta)^{\sigma/2}u, u_{t})(t,\cdot)\right\|_{L^{2}(\mathbb{R}^{n})} \lesssim (1+t)^{-\frac{n}{\sigma}\left(\frac{1}{m}-\frac{1}{2} \right)}\left\|u_{1}\right\|_{L^{m}(\mathbb{R}^{n})\cap L^{2}(\mathbb{R}^{n})}.
	\end{equation*}
\end{theorem}
\begin{remark}
	The conditions (\ref{11})-(\ref{12}) are assumed to get the same decay estimates of the semi-linear model with those of the corresponding linear model $(N(u)\equiv 0)$. More precisely, the bounds (\ref{11}) on $p+q$ and $n$ appear due to the application of Gagliardo-Nirenberg inequality from Lemma \ref{FGN} and the boundedness property of Riesz potential from Lemma \ref{hardy sobolev inequality}. The nonlocal (in space) nonlinearity does not cause any loss of decay rate in the semi linear model, as it does with the nonlocal (in time) nonlinearity known as nonlinear memory given by $$N(u)(t,x)=\int_{0}^{t}(t-s)^{-\gamma}|u(s,x)|^{p}ds, \ \ \gamma\in(0,1), \ \ p>1.$$
	In our forthcoming paper, we discuses the influence of modified nonlinear memory given by
	$$N(u)(t,x)=|u(t,x)|^{p}\int_{0}^{t}(t-s)^{-\gamma}|u(s,x)|^{q}ds, \ \ \gamma\in(0,1), \ \ p,q>1.$$
	on qualitative properties of solutions as in Theorem \ref{GlobalExistence1}.
	
\end{remark}
\begin{remark}
Because of the restriction imposed by Lemma \ref{hardy sobolev inequality}, we cannot deal with the Cauchy problem (\ref{1}) and nonlinearity (\ref{3}) if the initial velocity belongs to $L^{1}$ space, thus we choose it in $L^{m}$ space with $m\in(1,2)$. The H\"{o}lder inequality of the product of two functions can be used to solve this problem for generalized Hartree nonlinearity even if the initial velocity in $L^{1}$ space. This may be an important difference between the influence of (\ref{2}) and (\ref{3}).
\end{remark}
\begin{remark}
 We confirm that this result coincide with that in \cite{dabbicco hartree} for the \textit{same} parameters $p,q>1$ and $\alpha \in (0,n)$.
\end{remark}
\begin{proof} 
	
	Since we are dealing with semi-linear Cauchy problems, we use the Banach's fixed point theorem inspired from the book \cite[Page 303]{ebert reissig}.
	
	Here, we need to define a family of evolution spaces $X(T)$ for any $T>0$ with suitable norm $\|\cdot\|_{X(T)}$, also, we define an operator $O :u\in X(T)\longmapsto Ou(t,x)$ by $$ Ou(t,x)=G_{\sigma}(t,x)\ast u_{1}(x)+  \int_{0}^{t}G_{\sigma}(t-\tau,x)\ast I_{\alpha}(|u(\tau,x)|^{p+q})d\tau.$$
	If this operator satisfies the two inequalities: 
	\begin{align}
	\|Ou\|_{X(T)} &\lesssim \left\|u_{1}\right\|_{L^{m}(\mathbb{R}^{n})\cap L^{2}(\mathbb{R}^{n})}+\|u\|_{X(T)}^{p+q}, \ \  \forall u \in X(T), \label{13} \\ 
	\|Ou-Ov\|_{X(T)} &\lesssim \|u-v\|_{X(T)}\Big( \|u\|_{X(T)}^{p+q-1}+\|v\|_{X(T)}^{p+q-1}\Big), \ \  \forall u, v \in X(T), \label{14} 
	\end{align} 
	then, one can deduce the existence and uniqueness of a global (in time) solutions of (\ref{1})-(\ref{3}). Here, the smallness of the initial data $\left\|u_{1}\right\|_{L^{m}(\mathbb{R}^{n})\cap L^{2}(\mathbb{R}^{n})}<\varepsilon_0$ imply that the operator $O$ maps balls of $X(T)$ into balls of $X(T)$.
	
	The Banach space $X(T)$ is now defined for all $T>0$ as follows: 
	$$
	X(T):=\mathcal{C}\left([0,T], H^{\sigma}\right)\cap \mathcal{C}^{1}\left([0,T],L^{2}\right).
	$$
	
	We can equip this space with the norm 
	\begin{align*}
	\|u\|_{X(T)}&=\sup_{0\leq t\leq T}\Big((1+t)^{\frac{n}{\sigma}\left(\frac{1}{m}-\frac{1}{2} \right) -1}\|u(t,\cdot)\|_{L^{2}}\\
	&\hspace*{3cm}+(1+t)^{ \frac{n}{\sigma}\left(\frac{1}{m}-\frac{1}{2} \right)}\|((-\Delta)^{\sigma/2}u, u_{t})(t,\cdot)\|_{L^{2}}\Big).
	\end{align*} 
	
	By using the linear estimates from Lemma \ref{linear estimates}, it is clear the function $$u^{L}(t,x)=G_{\sigma}(t,x)\ast u_{1}(x)$$ that solves the linear model in (\ref{1}) belongs to $X(T)$ and we have
	\begin{align*}
	\|u^{L}\|_{X(T)}&=\sup_{0\leq t\leq T}\Big((1+t)^{\frac{n}{\sigma}\left(\frac{1}{m}-\frac{1}{2} \right)-1}\|u^{L}(t,\cdot)\|_{L^{2}}
	\\
	&\hspace*{3cm}+(1+t)^{ \frac{n}{\sigma}\left(\frac{1}{m}-\frac{1}{2} \right)}\|((-\Delta)^{\sigma/2}u^{L}, u_{t}^{L})(t,\cdot)\|_{L^{2}}\Big) \nonumber \\
	&
	\lesssim \left\|u_{1}\right\|_{L^{m}(\mathbb{R}^{n})\cap L^{2}(\mathbb{R}^{n})}.
	\end{align*} 
	To conclude inequality (\ref{13}), we prove
	\begin{equation}\label{15}
	\|u^{N}\|_{X(T)}\lesssim \|u\|_{X(T)}^{p+q}.
	\end{equation}
	We divide the interval $[0,t]$ into two sub-intervals $[0,t/2]$ and $[t/2,t]$ where we use the $L^{m}-L^{2}$ linear estimates if $\tau\in[0,t/2]$ and $L^{2}-L^{2}$ estimates if $\tau\in[t/2,t]$. From Lemma \ref{linear estimates} we estimate:
	\begin{align*}
	\|u^{N}(t,\cdot)\|_{L^{2}}&\lesssim\int_{0}^{t/2}(1+t-\tau)^{1-\frac{n}{\sigma}\left(\frac{1}{m}-\frac{1}{2} \right)} \left\|I_{\alpha}(|u|^{p+q})(\tau,\cdot)\right\| _{L^{m}\cap L^{2}} d\tau \nonumber \\
	&\hspace{3cm} 
	+\int_{t/2}^{t}(1+t-\tau)\left\|I_{\alpha}(|u|^{p+q})(\tau,\cdot)\right\|_{L^{2}}d\tau, 
	\end{align*}
	\begin{align*}
	\|((-\Delta)^{\sigma/2}u^{N},u_{t}^{N})(t,\cdot)\|_{L^{2}}	&\lesssim\int_{0}^{t/2}(1+t-\tau)^{-\frac{n}{\sigma}\left(\frac{1}{m}-\frac{1}{2} \right)} \left\|I_{\alpha}(|u|^{p+q})(\tau,\cdot)\right\| _{L^{m}\cap L^{2}}d\tau \nonumber \\
	&\hspace{2cm} 
	+\int_{t/2}^{t}\left\|I_{\alpha}(|u|^{p+q})(\tau,\cdot)\right\|_{L^{2}}d\tau.	
	\end{align*}
	Now, to control these integrals, we use the boundedness property of Riesz potential from Lemma \ref{hardy sobolev inequality} as follows:
	$$\|I_{\alpha}(|u|^{p+q})(\tau,\cdot)\|_{L^{m}}\lesssim \|u(\tau,\cdot)\|^{p+q}_{L^{q_{1}}}, \ \ \|I_{\alpha}(|u|^{p+q})(\tau,\cdot)\|_{L^{2}}\lesssim \|u(\tau,\cdot)\|^{p+q}_{L^{q_{2}}},$$
	where 
	$$q_{1}=\frac{mn(p+q)}{n+m\alpha}, \ \ \ \ q_{2}=\frac{2n(p+q)}{n+2\alpha}.$$ 
	Here, we know from the definition of the norm $\|u\|_{X(T)}$ that:
	\begin{equation*}
	(1+\tau)^{\frac{n}{\sigma}\left(\frac{1}{m}-\frac{1}{2}\right)}\|(-\Delta)^{\sigma/2}u(\tau,\cdot)\|_{L^{2}}\lesssim\|u\|_{X(T)},
	\end{equation*}
	\begin{equation*}
	(1+\tau)^{\frac{n}{\sigma}\left(\frac{1}{m}-\frac{1}{2}\right)-1}\|u(\tau,\cdot)\|_{L^{2}}\lesssim\|u\|_{X(T)}.
	\end{equation*}
	So, by the fractional Gagliardo-Nirenberg inequality from Lemma \ref{FGN}, we can estimate the above norms as follows:
	\begin{equation*}
	\left\|u(\tau,\cdot)\right\| _{L^{\frac{n(p+q)s}{n+s\alpha}}}^{p+q}\lesssim(1+\tau)^{-\Big((p+q)\left( \frac{n}{m\sigma}-1\right) -\frac{n}{s\sigma}-\frac{\alpha}{\sigma}\Big)}\|u\|_{X(T)}^{p+q}, \ s=m, 2,
	\end{equation*}
	provided that the conditions (\ref{11}) are satisfied for $p+q$ and $n$. Hence, we conclude 
	\begin{equation*}
	\left\|u(\tau,\cdot)\right\|^{p+q}_{L^{\frac{2n(p+q)}{n+2\alpha}}}+\left\|u(\tau,\cdot)\right\|^{p+q}_{L^{\frac{mn(p+q)}{n+m\alpha}}}\lesssim (1+\tau)^{-\Big((p+q)\left( \frac{n}{m\sigma}-1\right) -\frac{n}{m\sigma}-\frac{\alpha}{\sigma}\Big)}\|u\|_{X(T)}^{p+q}. 
	\end{equation*}
	By these equivalences
	$$(1+t-\tau)\approx (1+t) \ \text{if} \ \tau \in[0,t/2],\  \ (1+\tau)\approx (1+t) \ \text{if} \ \tau \in[t/2,t], $$
	and Lemma \ref{Integral inequality}, we estimates the first integral of $u^{N}$ over $[0,t/2]$ as follows
	$$
	\int_{0}^{t/2}(1+t-\tau)^{1-\frac{n}{\sigma}\left(\frac{1}{m}-\frac{1}{2}\right)}(1+\tau)^{-((p+q)\left( \frac{n}{m\sigma}-1\right) -\frac{n}{m\sigma}-\frac{\alpha}{\sigma})}\|u\|_{X(T)}^{p+q}d\tau$$
	$$
	\lesssim\|u\|_{X(T)}^{p+q}(1+t)^{1-\frac{n}{\sigma}\left(\frac{1}{m}-\frac{1}{2}\right)}\int_{0}^{t/2}(1+\tau)^{-((p+q)\left( \frac{n}{m\sigma}-1\right) -\frac{n}{m\sigma}-\frac{\alpha}{\sigma})}d\tau
	$$
	$$
	\lesssim\|u\|_{X(T)}^{p+q}(1+t)^{1
		-\frac{n}{\sigma}\left(\frac{1}{m}-\frac{1}{2}\right)},
	$$
	provided that $p+q>1+\frac{2m\sigma+m\alpha}{n-m\sigma}$. For the second integral over $[t/2,t]$ we have:
	
	$$\int_{t/2}^{t}(1+t-\tau)(1+\tau)^{-((p+q)\left( \frac{n}{m\sigma}-1\right) -\frac{n}{2\sigma}-\frac{\alpha}{\sigma})}\|u\|_{X(T)}^{p+q}d\tau$$
	$$\lesssim(1+t)^{2-((p+q)\left( \frac{n}{m\sigma}-1\right) -\frac{n}{2\sigma}-\frac{\alpha}{\sigma})}\|u\|_{X(T)}^{p+q}$$
	$$\lesssim(1+t)^{1-\frac{n}{\sigma}\left(\frac{1}{m}-\frac{1}{2}\right)}\|u\|_{X(T)}^{p+q},$$
	thanks to (\ref{12}), we arrive to the desired estimate for $u^{N}$ 
	$$
	(1+t)^{\frac{n}{\sigma}\left(\frac{1}{m}-\frac{1}{2}\right)-1}\|u^{N}(t,\cdot)\|_{L^{2}}\lesssim\|u\|_{X(T)}^{p+q}.
	$$
	Similarly, we can prove again
	$$ 
	(1+t)^{\frac{n}{\sigma}\left(\frac{1}{m}-\frac{1}{2}\right)}\|((-\Delta)^{\sigma/2}u^{N}, u_{t}^{N})(t,\cdot)\|_{L^{2}}\lesssim\|u\|_{X(T)}^{p+q}.
	$$
	Inequality (\ref{15}) is now proved, which implies (\ref{13}).
	
	To prove (\ref{14}) we choose two elements $u$, $v$ belong to $X(T)$, and we write
	$$Ou-Ov=\int_{0}^{t}G_{\sigma}(t-\tau,x)\ast I_{\alpha}(|u(\tau,x)|^{p+q}-|v(\tau,x)|^{p+q})d\tau.$$
	We divide this integral as above, thanks to the linearity and the boundedness property of Reisz potential, we have 
	$$\|I_{\alpha}(|u(\tau,x)|^{p+q}-|v(\tau,x)|^{p+q})\|_{L^{2}}\lesssim \||u(\tau,x)|^{p+q}-|v(\tau,x)|^{p+q}\|_{L^{\frac{2n}{n+2\alpha}}},$$
	$$ \|I_{\alpha}(|u(\tau,x)|^{p+q}-|v(\tau,x)|^{p+q})\|_{L^{m}}\lesssim \||u(\tau,x)|^{p+q}-|v(\tau,x)|^{p+q}\|_{L^{\frac{mn}{n+m\alpha}}}.$$ 	Now, by the H\"{o}lder's inequality, we derive  for $r=\frac{mn}{n+m\alpha}, \frac{2n}{n+2\alpha}$ the following
	\begin{equation*}
	\| |u(\tau,\cdot)|^{p+q}-|v(\tau,\cdot)|^{p+q}\|_{L^{r}}\lesssim \|u(\tau,\cdot)-v(\tau,\cdot)\|_{L^{r(p+q)}}\left(\|u(\tau,\cdot)\|_{L^{r(p+q)}}^{p+q-1}+\|v(\tau,\cdot)\|_{L^{r(p+q)}}^{p+q-1} \right).
	\end{equation*}
	Using the definition of the norm $\|u-v\|_{X(T)}$ and fractional Gagliardo-Nirenberg inequality we can prove without difficulty the inequality (\ref{14}). Hence, Theorem \ref{GlobalExistence1} is proved.	
\end{proof}
\section{Conclusion}\label{conclusion}
In this study, we proved that both nonlinearities (\ref{2}) and (\ref{3}) have distinct forms but give the same semi-critical exponent for global existence, where the first one is integrable and the second is not. It has been shown that the Hartree type nonlinearity appears naturally in various interesting physical systems such as  Schr\"{o}dinger equation of quantum mechanics coupled with Newton’s gravitational law or Penrose system (see [3] in \cite{dabbicco hartree}). If we use the second nonlinearity in original Penrose system it becomes
\begin{equation}
\left\lbrace 
\begin{array}{lll}
\hslash^{2}\Delta \psi(x)-V(x) \psi(x)+U(\psi)(x)=0\ \ \ \hfill& x\in \mathbb{R}^{3},
&\cr
\\
\hslash^{2}\Delta U(x)+|\psi(x)|^{3}=0\ \ \ \hfill& x\in \mathbb{R}^{3},
\end{array}\right.  
\end{equation}
This system can be explicitly solved with respect to $U$, so that it becomes the following single nonlinear nonlocal stationary equation
$$\hslash^{2}\Delta \psi(x)-V(x) \psi(x)+\frac{1}{4\pi\hslash^{2}}\Big(\int_{\mathbb{R}^{3}}\frac{|\psi(y)|^{3}}{|x-y|}dy\Big)=0, \ \ \ \ x\in \mathbb{R}^{3}.$$
The original Penrose equation is given by (see [3] in \cite{dabbicco hartree})
$$\hslash^{2}\Delta \psi(x)-V(x) \psi(x)+\frac{1}{4\pi\hslash^{2}}\Big(\int_{\mathbb{R}^{3}}\frac{|\psi(y)|^{2}}{|x-y|}dy\Big)\psi(x)=0, \ \ \ \ x\in \mathbb{R}^{3}.$$
The same idea can be used in other physical models.	


\begin{thebibliography}{999999}
	\bibitem{dabbicco hartree} D'Abbicco. M, Global Small Data Solutions for an
	Evolution Equation with Structural
	Damping and Hartree-Type Nonlinearity. In: Kähler, U., Reissig, M., Sabadini, I., Vindas, J. (eds) Analysis, Applications, and Computations. ISAAC 2021. Trends in Mathematics(). Birkhäuser, Cham. https://doi.org/10.1007/978-3-031-36375-7$\_$43 
	\\
\bibitem{pham reissig kainane} Pham. D. T, Kainane Mozadek. M, Reissig. M, Global existence for semi-linear structurally damped $\sigma$-evolution models, J. Math. Anal. Appl, 431(2015), 569--596.
\\
   \bibitem{ebert reissig} Ebert. M. R, Reissig. M, Methods for partial differential equations, qualitative properties of solutions, phase space analysis, semilinear models. Birkh\"{a}user, 2018.
   \\
   \bibitem{stein} Stein. E. M, Singular Integrals and Differentiability Properties of Functions, (Princeton Math. Ser. 30, Princeton, N.J., 1970).


\end{thebibliography}
\end{document}